\pdfoutput=1
\documentclass{article}
\usepackage[english]{babel}
\usepackage{amsmath,amssymb,amsthm}
\usepackage{amscd}
\usepackage{graphicx}
\usepackage{braket}
\usepackage{tikz-cd}
\usepackage[hidelinks]{hyperref}

\providecommand{\longhookrightarrow}{\lhook\joinrel\longrightarrow}

\theoremstyle{plain}
\newtheorem{theorem}{Theorem}[section]
\newtheorem{lemma}[theorem]{Lemma}
\newtheorem{proposition}[theorem]{Proposition}
\newtheorem{corollary}[theorem]{Corollary}
\theoremstyle{definition}
\newtheorem{definition}[theorem]{Definition}
\newtheorem{notation}[theorem]{Notation}
\theoremstyle{remark}
\newtheorem*{remark}{Remark}
\newtheorem*{example}{Example}

\title{Ultraproducts of crossed product C$^*$algebras}
\author{Zhengyu Fu\thanks{This is the author's accepted manuscript of an article published in
Canad.\ Math.\ Bull.\ (2026), \url{https://doi.org/10.4153/S0008439526102070}.}}
\date{}

\begin{document}
\maketitle

\begin{abstract}
We study the relationship between the ultraproduct of a crossed product C$^*$algebra
$(\mathcal{A}\rtimes_{r}G)^{\omega}$ and the crossed product of an ultraproduct C$^*$algebra
$\mathcal{A}^{\omega}\rtimes_{r}G$ for a fixed free ultrafilter $\omega$ on $\mathbb{N}$.
For a general locally compact amenable group $G$, we show that
$\mathcal{A}^{\omega}_{\alpha}\rtimes_{r}G$ embeds into $(\mathcal{A}\rtimes_{r}G)^{\omega}$.
This problem was first introduced by Reiji Tomatsu in the framework of von Neumann algebras
in \cite{Tom17}.
\end{abstract}



\section{Introduction}

Ultraproducts play an important role in the study of operator algebras and their
automorphisms (cf.\ \cite{AH14}). In this context, Reiji Tomatsu initiated the study of
Rohlin flows on von Neumann algebras using ultraproduct techniques \cite{MT16} and recognized
a close relationship between ultraproducts and crossed products, leading to the following
theorem.

\begin{theorem}[{\cite[Theorem 3.6]{Tom17}}]\label{thm:Tomatsu}
Let $\alpha$ be an action of a second countable locally compact Hausdorff group $G$ on a
von Neumann algebra $M$ with separable predual. Then the following statements hold:
\begin{itemize}
\item[(1)] There exists a canonical embedding from
$M^{\omega}_{\alpha}\rtimes_{\alpha^{\omega}}G$ into $(M\rtimes_{\alpha}G)^{\omega}$.
\item[(2)] If $G$ is abelian, then the embedding induces an isomorphism from
$M^{\omega}_{\alpha}\rtimes_{\alpha^{\omega}}G$ onto $(M\rtimes_{\alpha}G)^{\omega}_{\hat{\alpha}}$.
\end{itemize}
\end{theorem}

Theorem~\ref{thm:Tomatsu} provides a powerful tool for analyzing the structure of ultraproduct
von Neumann algebras through continuous or discrete crossed product decompositions. As an
application, let $M=N\rtimes_{\theta}\mathbb{R}$ be the continuous crossed product
decomposition of a properly infinite von Neumann algebra $M$, where $N$ is a semifinite von
Neumann algebra and $\theta$ is a trace scaling one-parameter automorphism group. Then
Theorem~\ref{thm:Tomatsu} implies that the ultraproduct von Neumann algebra $M^{\omega}$ also
admits a continuous crossed product decomposition,
$M^{\omega}=N^{\omega}\rtimes_{\theta^{\omega}}\mathbb{R}$. Hence, the continuous
decomposition structure is preserved under ultraproducts, allowing one to lift the analysis of
the flow of weights to ultraproduct von Neumann algebras.

The main difficulty is that continuous group actions do not continuously extend to
ultraproduct C$^*$algebras. This leads us to the notion of continuous part introduced by
Kishimoto \cite{Kis96}. For a C$^*$algebra $\mathcal{A}$, we simply write
$\ell^{\infty}(\mathcal{A})=\ell^{\infty}(\mathbb{N},\mathcal{A})$. Define
$C^{\omega}(\mathcal{A}) = \{x \in \ell^{\infty}(\mathcal{A}) : \lim_{n \to \omega} \|x(n)\| = 0 \}
\triangleleft\ell^{\infty}(\mathcal{A})$ and the ultraproduct of $\mathcal{A}$ is given by
$\mathcal{A}^{\omega} = \ell^{\infty}(\mathcal{A}) / C^{\omega}(\mathcal{A})$. Given a continuous
action $\alpha: G\curvearrowright\mathcal{A}$, it naturally induces an algebraic ultra-action
$\alpha^{\omega}: G \to\mathrm{Aut}(\mathcal{A}^{\omega})$ on the ultrapower (note that
$\mathcal{A}^{\omega}$ is non-separable in general and $\mathrm{Aut}(\mathcal{A}^{\omega})$ may
not be a Polish group). By restricting to the continuous part
$\mathcal{A}^{\omega}_{\alpha}:=\ell^{\infty}_{\alpha}(\mathcal{A})/C^{\omega}(\mathcal{A})$, we
then obtain a C$^*$ dynamical system. Here,
$\ell^{\infty}_{\alpha}(\mathcal{A}):=\{x\in\ell^{\infty}(\mathcal{A}):
s\mapsto \bar{\alpha}_{s}(x) \text{ is continuous}\}$, where $\bar{\alpha}$ acts on
$\ell^{\infty}(\mathcal{A})$ in the natural way. In this paper, we first consider the case of a
discrete group, which we always denote by $\Gamma$. We conclude with the following theorems.

\begin{theorem}[Theorems~\ref{thm:embedding} and~\ref{thm:reverse}]\label{thm:main}
Let $\alpha$ be an action of an exact discrete group $\Gamma$ on a separable  C$^*$algebra $\mathcal{A}$.
Then the following statements hold:
\begin{itemize}
\item[(1)] There exists a canonical embedding from
$\mathcal{A}^{\omega}_{\alpha}\rtimes_{\alpha^{\omega},r}\Gamma$ into
$(\mathcal{A}\rtimes_{\alpha,r}\Gamma)^{\omega}$.
\item[(2)] If $\Gamma$ is abelian, then the embedding induces an isomorphism from
$\mathcal{A}^{\omega}_{\alpha}\rtimes_{\alpha^{\omega}}\Gamma$ onto
$(\mathcal{A}\rtimes_{\alpha}\Gamma)^{\omega}_{\hat{\alpha}}$.
\end{itemize}
\end{theorem}

\section{Preliminaries}

\begin{definition}
Let $\Gamma$ be a discrete group and let $\lambda$ be the left regular representation of
$\Gamma$ on $\ell^{2}(\Gamma)$. For a $\Gamma$-C$^*$algebra
$\mathcal{A}\stackrel{\pi}{\subset}B(\mathcal{H})$, define the representation $\tilde{\pi}$ of
$\mathcal{A}$ on $\mathcal{H}\otimes\ell^{2}(\Gamma)$ by
$\tilde{\pi}(a)(\xi\otimes\delta_{g})=(\alpha_{g^{-1}}(a)\xi)\otimes\delta_{g}$. We can see that
$(\tilde{\pi},1\otimes\lambda,\mathcal{H}\otimes \ell^{2}(\Gamma))$ is a covariant
representation. We define the reduced crossed product $\mathcal{A}\rtimes_{\alpha,r}\Gamma$ by
$\overline{\mathrm{span}}\{\tilde{\pi}(a)(1\otimes\lambda_{s}):a\in\mathcal{A}, s\in \Gamma\}
\subset B(\mathcal{H}\otimes \ell^{2}(\Gamma))$.
\end{definition}

\begin{definition}\label{def:exact}
A C$^*$algebra $\mathcal{A}$ is said to be exact if for every short exact sequence of
C$^*$-algebras $0\to I\to E\to Q\to0$, the sequence
$0\to I\otimes_{\mathrm{min}}\mathcal{A}\to E\otimes_{\mathrm{min}}\mathcal{A}\to
Q\otimes_{\mathrm{min}}\mathcal{A} \to 0$ is exact. Moreover, a discrete group is exact if its
reduced group C$^*$algebra is exact.
\end{definition}

For an exact group, the reduced crossed product preserves short exact sequences (see also
\cite[Theorem 5.1.10]{BO08}). Consequently, for the C$^*$short exact sequence, the sequence
\[
0\to C^{\omega}\rtimes_{r}\Gamma\stackrel{\iota\rtimes\mathrm{id}}{\longrightarrow}
\ell^{\infty}(\mathcal{A})\rtimes_{r}\Gamma\stackrel{q\rtimes\mathrm{id}}{\longrightarrow}
\ell^{\infty}(\mathcal{A})/C^{\omega}\rtimes_{r}\Gamma\to0
\]
is exact and hence induces a well-defined $*$isomorphism
\[
\frac{\ell^{\infty}(\mathcal{A})\rtimes_{r}\Gamma}{C^{\omega}\rtimes_{r}\Gamma}
\stackrel{\sim}{\longrightarrow}\ell^{\infty}(\mathcal{A})/C^{\omega}\rtimes_{r}\Gamma.
\]

\subsection{Crossed product for locally compact group}

(See \cite[Chapter 2]{DW07}.) Let $G$ be a locally compact Hausdorff group, which is always
assumed to be second countable. Suppose that $(\mathcal{A},G,\alpha)$ is a C$^*$ dynamical
system and that, for each $f\in C_{c}(G,\mathcal{A})$, we define
\[
\|f\|_{u}:=\sup\{\|\pi\rtimes U(f)\|:(\pi,U)\in \mathrm{Covrep}(\mathcal{A},G)\},
\]
where $(\pi,U)$ is a covariant representation of $(\mathcal{A},G)$ and $\pi\rtimes U$ is the
integrated form of $(\pi,U)$. The crossed product $\mathcal{A}\rtimes_{\alpha}G$ is the
completion of $C_{c}(G,\mathcal{A})$ with respect to $\|\cdot\|_{u}$. As is well known, there is
a one-to-one correspondence between nondegenerate covariant representations of
$(\mathcal{A},G,\alpha)$ and nondegenerate representations of $\mathcal{A}\rtimes_{\alpha}G$.
Another well-known fact is that, if $G$ is amenable, then $\mathcal{A}\rtimes_{\alpha}G$
coincides with the reduced crossed product $\mathcal{A}\rtimes_{\alpha,r}G$.

\subsection{Duality}

(See \cite[Theorem 6]{Rae88}.) For a locally compact abelian group $G$,
\[
\hat{\alpha}:\hat{G}\to\mathrm{Aut}(C_{c}(G,\mathcal{A}));\quad
\hat{\alpha}_{\gamma}(x)(s):=\braket{s,\gamma}x(s)
\]
extends to a continuous action $\hat{\alpha}:\hat{G}\curvearrowright\mathcal{A}\rtimes_{\alpha}G$.
By the Takai duality, we have an isomorphism from
$(\mathcal{A}\rtimes_{\alpha}G)\rtimes_{\hat{\alpha}}\hat{G}$ onto
$\mathcal{A}\otimes \mathbb{K}(L^{2}(G))$ such that $\hat{\hat{\alpha}}$ is carried into
$\alpha\otimes \mathrm{Ad}\rho$, where $\rho$ denotes the right regular representation.

For C$^*$algebra-valued functions, one can avoid issues with approximation by simple functions
by using a more direct approach than the Bochner integral. The following well-known lemma
provides a convenient definition of integration and its basic properties.

\begin{lemma}[{See \cite[Lemma 7]{Rae88}}]\label{lem:integral}
For $f\in C_{c}(G,\mathcal{A})$, there exists a unique element $\int f(s)\,ds\in\mathcal{A}$ such
that for any nondegenerate representation $\pi$ of $\mathcal{A}$,
\[
\Braket{\pi\Bigl(\int f(s)\,ds\Bigr)\xi,\eta}=\int\braket{\pi(f(s))\xi,\eta}\,ds
\quad\text{for } \xi,\eta\in \mathcal{H}_{\pi},
\]
where $ds$ denotes the Haar measure on $G$. Further, the following properties hold:
\begin{itemize}
\item $(\int f(s)\,ds)\,a = \int f(s)a\,ds$ for $a\in\mathcal{A}$;
\item $(\int f(s)\,ds)^{*} = \int f(s)^{*}\,ds$;
\item $\varphi(\int f(s)\,ds) = \int \varphi(f(s))\,ds$ for any homomorphism
$\varphi:\mathcal{A}\to\mathcal{B}$.
\end{itemize}
\end{lemma}

\section{Main result}

\subsection{On embedding relations}

\begin{lemma}\label{lem:linfty}
Let $\Gamma$ be a discrete group. Then $\ell^{\infty}(\mathcal{A})\rtimes_{\bar{\alpha},r}\Gamma$ is
a C$^*$subalgebra of $\ell^{\infty}(\mathcal{A}\rtimes_{\alpha,r}\Gamma)$. (An action
$\bar{\alpha}$ on $\ell^{\infty}(\mathcal{A})$ is naturally considered as
$\bar{\alpha}_{g}((x_{n})_{n})=(\alpha_{g}(x_{n}))_{n}$.)
\end{lemma}

\begin{proof}
The map
\[
C_{c}(\Gamma,\ell^{\infty}(\mathcal{A}))\ni\sum_{s\in F;\,F\Subset\Gamma}(a_{s}(n))_{n}u_{s}
\mapsto \Bigl(\sum_{s\in F}a_{s}(n)u_{s}\Bigr)_{n}\in \ell^{\infty}(C_{c}(\Gamma,\mathcal{A}))
\]
is a well-defined $*$homomorphism. Let $\iota:\mathcal{A}\hookrightarrow B(\mathcal{H})$ be a
faithful representation. Then the $\ell^{\infty}$-direct sum of $\iota$, denoted by
$\bigoplus_{n}\iota$, defines a faithful representation of $\ell^{\infty}(\mathcal{A})$.
Consequently, the induced representation
\[
\widetilde{\bigoplus_{n}\iota}:\ell^{\infty}(\mathcal{A}) \hookrightarrow
B(\mathcal{H}^{\oplus_{n}}\otimes \ell^{2}(\Gamma))
\]
is also faithful. Since
$\mathcal{H}^{\oplus_{n}}\otimes \ell^{2}(\Gamma)\simeq\bigoplus_{n}(\mathcal{H}\otimes \ell^{2}(\Gamma))$
via a unitary isomorphism, it is straightforward to verify that
$\widetilde{\bigoplus_{n}\iota}\rtimes\oplus(1\otimes\lambda)$ and
$\bigoplus_{n}(\tilde{\iota}\rtimes(1\otimes\lambda))$ are unitarily equivalent. Thus, we compute
\begin{align*}
\sup_{n}\Bigl\|\sum_{s\in F}a_{s}(n)u_{s}\Bigr\|_{\mathrm{min}}
&=\sup_{n}\Bigl\|\tilde{\iota}\rtimes(1\otimes\lambda)\sum_{s}a_{s}(n)u_{s}\Bigr\|\\
&=\Bigl\|\widetilde{\bigoplus_{n}\iota}\rtimes\oplus(1\otimes\lambda)\sum_{s}(a_{s}(n))_{n}u_{s}\Bigr\|.
\end{align*}
Therefore, the $*$homomorphism defined above extends to an injective homomorphism from
$\ell^{\infty}(\mathcal{A})\rtimes_{\bar{\alpha},r}\Gamma$ into
$\ell^{\infty}(\mathcal{A}\rtimes_{\alpha,r}\Gamma)$.
\end{proof}

Let us denote $C^{\omega}(\mathcal{A}\rtimes_{\alpha,r}\Gamma)$ by $\bar{C}^{\omega}$. For
$a\in C^{\omega}$, $\lim_{n\to\omega}\|\bar{\alpha}_{g}(a)(n)\|\leq\lim_{n\to\omega}\|a(n)\|$
$({}^{\forall}g\in \Gamma)$ implies $\bar{\alpha}_{g}(C^{\omega})=C^{\omega}$
$({}^{\forall}g\in \Gamma)$, and thus $\bar{\alpha}$ is also an action on $C^{\omega}$.
Therefore, the following relationship exists:
\[
\begin{array}{ccc}
\ell^{\infty}(\mathcal{A})\rtimes_{\bar{\alpha},r}\Gamma &\hookrightarrow&
\ell^{\infty}(\mathcal{A}\rtimes_{\alpha,r}\Gamma)\\
\rotatebox{90}{$\subset$}& & \rotatebox{90}{$\subset$}\\
C^{\omega}\rtimes_{\bar{\alpha},r}\Gamma & \hookrightarrow & \bar{C}^{\omega}
\end{array}
\]

\begin{theorem}\label{thm:embedding}
Let $\Gamma$ be an exact discrete group. For any action $\alpha:\Gamma\curvearrowright\mathcal{A}$,
we have $\mathcal{A}^{\omega}\rtimes_{\alpha^{\omega},r}\Gamma\subset
(\mathcal{A}\rtimes_{\alpha,r}\Gamma)^{\omega}$.
\end{theorem}

\begin{proof}
Let $F\Subset\Gamma$. By Lemma~\ref{lem:linfty}, the map
\[
\mathcal{J}:\ell^{\infty}(\mathcal{A})\rtimes_{\bar{\alpha},r}\Gamma\ni
\sum_{s\in F}(a_{s}(n))_{n} u_{s}\mapsto \Bigl(\sum_{s\in F}a_{s}(n) u_{s}\Bigr)_{n}\in
\ell^{\infty}(\mathcal{A}\rtimes_{\alpha,r}\Gamma)
\]
is a $*$isometry. It is well known that, for an exact group, the reduced crossed product
preserves short exact sequences \cite[Chapter 5]{BO08}. Therefore, we obtain the following
double exact sequences:
\[
  \scalebox{0.92}{$\begin{CD}
     0 @>>>  C^{\omega}\rtimes_{\bar{\alpha},r}\Gamma @>>>  \ell^{\infty}(\mathcal{A})\rtimes_{\bar{\alpha},r}\Gamma  @>{q\rtimes\mathrm{id}}>>  \ell^{\infty}(\mathcal{A})/C^{\omega}\rtimes_{\alpha^{\omega},r}\Gamma  @>>>  0 \\
    @.     @.    @VV{\mathcal{J}}V  @VV{\Theta}V   @. \\
     0 @>>> \bar{C}^{\omega} @>>>  \ell^{\infty}(\mathcal{A}\rtimes_{\alpha,r}\Gamma) @>\bar{q}>>  \ell^{\infty}(\mathcal{A}\rtimes_{\alpha,r}\Gamma)/\bar{C}^{\omega} @>>>  0
  \end{CD}$}
\]
Here, $q$ and $\bar{q}$ are the canonical quotient maps. As mentioned after
Definition~\ref{def:exact}, we note that
$\ell^{\infty}(\mathcal{A})/C^{\omega}\rtimes_{\alpha^{\omega},r}\Gamma$ is nothing but
$\ell^{\infty}(\mathcal{A})\rtimes_{\bar{\alpha},r}\Gamma/C^{\omega}\rtimes_{\bar{\alpha},r}\Gamma$.
We define a $*$homomorphism $\Theta$ so that the right-hand side of the above diagram commutes,
that is,
\[
\Theta\circ q\rtimes\mathrm{id}:=\bar{q}\circ\mathcal{J}.
\]
Furthermore, by Lemma~\ref{lem:linfty}, $C^{\omega}\rtimes_{\bar{\alpha},r}\Gamma\subset\bar{C}^{\omega}$
as a C$^*$subalgebra. Therefore, $\Theta$ is well-defined. The following argument is inspired by
a technique appearing in \cite{Pis09}. It remains to show that $\Theta$ is injective. Assume
$x\in \ell^{\infty}(\mathcal{A})/C^{\omega}\rtimes_{\alpha^{\omega},r}\Gamma$ and $\Theta(x)=0$. For
$y\in \ell^{\infty}(\mathcal{A})\rtimes_{\bar{\alpha},r}\Gamma$ such that $q\rtimes\mathrm{id}(y)=x$,
we have $\mathcal{J}(y)\in\bar{C}^{\omega}$, i.e.,
\[
{}^{\forall}\varepsilon>0,\ {}^{\exists}S\in\omega \text{ s.t. }
\sup_{n\in S}\|\mathcal{J}(y)_{n}\|_{\mathcal{A}\rtimes_{\alpha,r}\Gamma}<\varepsilon.
\]
We define $p_{S}:\ell^{\infty}(\mathcal{A})\to \ell^{\infty}(\mathcal{A})$ as follows: for
$a\in \ell^{\infty}(\mathcal{A})$, $p_{S}(a)_{n}=0$ if $n\notin S$, and the $n$-th component
remains the same if $n\in S$. Clearly, $p_{S}$ is a projection. For any
$(a_{n})_{n}\in \ell^{\infty}$, we have
\[
q(p_{S}((a_{n})_{n}))=q((a_{n})_{n}).
\]
Indeed, for a fixed $\varepsilon$ as above, by the definition of $S$, we have
$S\subset\{i:\|p_{S}(a_{i})-a_{i}\|<\varepsilon\}$, which implies that
$\{i:\|p_{S}(a_{i})-a_{i}\|<\varepsilon\} \in\omega$. Because $p_{S}$ is $\Gamma$-equivariant,
we can consider $p_{S}\rtimes \mathrm{id}_{\Gamma}$, which acts as the identity on
$\mathcal{A}^{\omega}\rtimes_{\alpha^{\omega},r}\Gamma$. Thus, estimating the norm of $x$, we have
\begin{align*}
\|x\|&=\|q\rtimes\mathrm{id}(y)\|=\|q\rtimes\mathrm{id}(p_{S}\rtimes \mathrm{id})(y)\|
\leq\|p_{S}\rtimes \mathrm{id}(y)\|_{\ell^{\infty}(\mathcal{A})\rtimes_{\bar{\alpha},r}\Gamma}\\
&=\|\mathcal{J}(p_{S}\rtimes \mathrm{id}(y))\|_{\ell^{\infty}(\mathcal{A}\rtimes_{\alpha,r}\Gamma)}
=\sup_{n\in S}\|\mathcal{J}(y)_{n}\|_{\mathcal{A}\rtimes_{r}\Gamma}<\varepsilon.
\end{align*}
Since $\varepsilon$ is arbitrary, we conclude $x=0$.
\end{proof}

\begin{remark}
If $\mathcal{B}$ is an exact C$^*$algebra, then $\mathcal{A}^{\omega}\otimes_{\mathrm{min}}\mathcal{B}$
is just the quotient C$^*$algebra of $\ell^{\infty}(\mathcal{A})\otimes_{\mathrm{min}}\mathcal{B}$, as
mentioned after Definition~\ref{def:exact}. Similarly, by the proof of
Theorem~\ref{thm:embedding}, we obtain the following corollary, which is well known to experts
and appears, at least implicitly, in \cite{BO08}.
\end{remark}

\begin{corollary}
If $\mathcal{B}$ is an exact C$^*$algebra, then
$\mathcal{A}^{\omega}\otimes_{\mathrm{min}}\mathcal{B}\subset(\mathcal{A}\otimes_{\mathrm{min}}\mathcal{B})^{\omega}$.
\end{corollary}

\begin{example}
The nuclearity of $\mathbb{K}(L^{2}(G))$ implies exactness, and so
\[
\mathcal{A}^{\omega}\otimes \mathbb{K}(L^{2}(G))\subset(\mathcal{A}\otimes \mathbb{K}(L^{2}(G)))^{\omega}
\]
holds.
\end{example}

\begin{proposition}\label{prop:finite}
Let $G$ be a finite group. Then, for any action $\alpha:G\curvearrowright\mathcal{A}$,
$\mathcal{A}^{\omega}\rtimes_{\alpha^{\omega},r}G=(\mathcal{A}\rtimes_{\alpha,r}G)^{\omega}$ holds.
\end{proposition}

\begin{proof}
By Lemma~\ref{lem:linfty}, it suffices to show that $\bar{C}^{\omega}\subset C^{\omega}\rtimes_{\bar{\alpha},r}G$.
Since $G$ is finite, every element of $\mathcal{A}\rtimes_{r}G$ can be written as a finite sum.
Hence, every element of $(\mathcal{A}\rtimes_{r} G)^{\omega}$ is represented by a bounded sequence
of such finite sums.

We can take the conditional expectation on the reduced crossed product,
$E:C^{\omega}\rtimes_{\bar{\alpha},r}G\to C^{\omega};\
\sum_{s}[(a_{s}(n))]_{\omega}u_{s}\mapsto [(a_{e}(n))]_{\omega}$, and
\[
\lim_{n\to\omega}\|a_{g}(n)\|=\lim_{\omega}\Bigl\|E\Bigl(u_{g}^{*}\sum_{s}(a_{s}(n))_{n}u_{s}\Bigr)\Bigr\|
\leq\lim_{\omega}\Bigl\|\sum_{s}a_{s}(n)u_{s}\Bigr\|_{\mathcal{A}\rtimes_{\bar{\alpha},r}G}.
\]
This completes the proof.
\end{proof}

In the case of a locally compact group, Lemma~\ref{lem:linfty} still holds. Indeed, there is a
well-defined $G$-equivariant homomorphism
\[
\Phi:C_{c}(G,\ell^{\infty}_{\alpha}(\mathcal{A}))\to \ell^{\infty}(C_{c}(G,\mathcal{A}));\quad
\Phi(f)(n)(s):=f(s)(n).
\]
However, ultralimits do not satisfy the dominated convergence theorem, and hence $\Theta$ cannot
be constructed as in Theorem~\ref{thm:embedding}. This is due to the fact that the norm estimate
required for the well-definedness of $\Theta$ involves an integral, and in general one cannot
interchange the ultralimit and the integral.

\begin{proposition}\label{prop:hom}
Let $G$ be a locally compact group. Then there exists a homomorphism from
$\mathcal{A}^{\omega}_{\alpha}\rtimes_{\alpha^{\omega}}G$ into $(\mathcal{A}\rtimes_{\alpha}G)^{\omega}$.
\end{proposition}

\begin{proof}
Let $(\iota,u)$ be the canonical covariant representation of $(\mathcal{A},G,\alpha)$ into the
multiplier algebra $M(\mathcal{A}\rtimes_{\alpha}G)$. Thanks to the natural inclusion
$M(\mathcal{A}\rtimes_{\alpha}G)^{\omega}\subset M((\mathcal{A}\rtimes_{\alpha}G)^{\omega})$, the
(restriction of the) ultraproduct representation
\[
\iota^{\omega}:\mathcal{A}^{\omega}_{\alpha}\hookrightarrow M((\mathcal{A}\rtimes_{\alpha}G)^{\omega}),\quad
u^{\omega}:G\hookrightarrow M((\mathcal{A}\rtimes_{\alpha} G)^{\omega}),
\]
where $u^{\omega}(g)=[u_{g}]_{\omega}$, is also covariant. However, the strict continuity of
$u^{\omega}$ may fail in general (see the Remark below). Define the strictly continuous part by
\[
\mathcal{U}:=\{x\in (\mathcal{A}\rtimes G)^{\omega}:s\mapsto u^{\omega}_{s}x \text{ and }
xu^{\omega}_{s} \text{ are continuous}\},
\]
which is a hereditary C$^*$subalgebra of $(\mathcal{A}\rtimes G)^{\omega}$.

We now show that $\iota^{\omega}(\mathcal{A}^{\omega}_{\alpha})\subset M(\mathcal{U})$. For
$x\in\mathcal{A}^{\omega}_{\alpha}$, by covariance, we have
$u^{\omega}_{s}\iota^{\omega}(x)=\iota^{\omega}(\alpha^{\omega}_{s}(x))u^{\omega}_{s}$, and thus,
for $y\in\mathcal{U}$,
\[
\|u^{\omega}_{s}\iota^{\omega}(x)y-\iota^{\omega}(x)y\|\leq
\|\iota^{\omega}(\alpha^{\omega}_{s}(x)-x)\|\|u^{\omega}_{s}y\|
+\|\iota^{\omega}(x)\|\|u^{\omega}_{s}y-y\|\to 0\quad(s\to e).
\]
Taking adjoints, we conclude that $\iota^{\omega}(\mathcal{A}^{\omega}_{\alpha})\subset M(\mathcal{U})$.
Similarly, by definition of $\mathcal{U}$, we have $u^{\omega}_{G}\subset M(\mathcal{U})$. Hence, by
universality, the integrated form $\iota^{\omega}\rtimes u^{\omega}$ yields a homomorphism from
$\mathcal{A}^{\omega}_{\alpha}\rtimes_{\alpha^{\omega}}G$ into $M(\mathcal{U})$. Actually, the image
of $\iota^{\omega}\rtimes u^{\omega}$ lies in $\mathcal{U}$, i.e., for
$a\in\mathcal{A}^{\omega}_{\alpha}$ and $f\in C_{c}(G)$,
$\iota^{\omega}\rtimes u^{\omega}(af)=\iota^{\omega}(a)u^{\omega}(f)\in\mathcal{U}$. This follows
from the observation that the maps
\begin{align*}
s&\mapsto u^{\omega}_{s}\iota^{\omega}(a)u^{\omega}(f)=\iota^{\omega}(\alpha^{\omega}_{s}(a))u^{\omega}(s.f),\\
s&\mapsto \iota^{\omega}(a)u^{\omega}(f)u^{\omega}_{s}=\iota^{\omega}(a)u^{\omega}(f.s)
\end{align*}
are norm continuous, where $s.f$ and $f.s$ denote the left and right translations,
respectively. Thus, we obtain the desired map from $\mathcal{A}^{\omega}_{\alpha}\rtimes G$ to
$(\mathcal{A}\rtimes G)^{\omega}$.
\end{proof}

\begin{remark}
Let $\lambda:\mathbb{R}\curvearrowright L^{2}(\mathbb{R})$ be the left regular representation and
let $\lambda^{\omega}:\mathbb{R}\to M(\mathbb{K}(L^{2}(\mathbb{R}))^{\omega})$ denote its
ultraproduct representation. Choose $\xi_{n}\in L^{2}(\mathbb{R})$ such that
$\mathrm{supp}\,\xi_{n}\subset[-\frac{1}{n},\frac{1}{n}]$ and $\|\xi_{n}\|_{2}=1$ for each
$n\in\mathbb{N}$. Then, for every fixed $t\in\mathbb{R}\setminus\{0\}$, we have
$\|\lambda_{t}\xi_{n}-\xi_{n}\|^{2}_{2}=2$ for all sufficiently large $n$. Define
$a:=[(\xi_{n}\otimes\xi_{n})]_{\omega}\in\mathbb{K}(L^{2}(\mathbb{R}))^{\omega}$. We obtain
$\|\lambda_{t}^{\omega}a-a\|=\lim_{n\to\omega}\|\lambda_{t}\xi_{n}-\xi_{n}\|_{2}=\sqrt{2}$, which
implies $\lim_{t\to 0}\lambda^{\omega}_{t}a\neq a$.
\end{remark}

For the reduced crossed product, we can faithfully represent
\begin{gather*}
\mathcal{A}^{\omega}_{\alpha}\rtimes_{\alpha^{\omega},r}G
\stackrel{\widetilde{\pi^{\omega}}\rtimes(1\otimes\lambda)}{\longhookrightarrow}
B(\mathcal{H}^{\omega}\otimes L^{2}(G))
\quad\text{and}\\
(\mathcal{A}\rtimes_{\alpha,r}G)^{\omega}
\stackrel{(\tilde{\pi}\rtimes(1\otimes\lambda))^{\omega}}{\longhookrightarrow}
B((\mathcal{H}\otimes L^{2}(G))^{\omega}).
\end{gather*}
Hence, it is crucial to consider the transformation
between the two Hilbert spaces, $\mathcal{H}^{\omega}\otimes L^{2}(G)$ and
$(\mathcal{H}\otimes L^{2}(G))^{\omega}$.

\begin{lemma}\label{lem:V}
The map $V:\mathcal{H}^{\omega}\otimes_{\mathrm{alg}}L^{2}(G)\to(\mathcal{H}\otimes L^{2}(G))^{\omega}$
defined by $V((\xi_{n})_{\omega}\otimes\eta)=(\xi_{n}\otimes\eta)_{\omega}$ extends to an isometry on
$\mathcal{H}^{\omega}\otimes L^{2}(G)$. Moreover, for each $a\in\mathcal{A}_{\alpha}^{\omega}$ and
$g\in G$, the following relations hold:
\[
V\widetilde{\pi^{\omega}}(a)=\tilde{\pi}^{\omega}(a)V,\qquad
V(1\otimes\lambda_{g})=(1\otimes\lambda_{g}^{\omega})V.
\]
\end{lemma}

\begin{proof}
This follows from a straightforward calculation.
\end{proof}

\begin{theorem}\label{thm:amenable}
If $G$ is amenable, then
$\mathcal{A}^{\omega}_{\alpha}\rtimes_{\alpha^{\omega}}G\subset(\mathcal{A}\rtimes_{\alpha}G)^{\omega}$.
\end{theorem}

\begin{proof}
Because $G$ is amenable, the above discussion allows us to faithfully represent
$\mathcal{A}^{\omega}_{\alpha}\rtimes_{\alpha^{\omega}}G$ on $\mathcal{H}^{\omega}\otimes L^{2}(G)$
by $\widetilde{\pi^{\omega}}\rtimes(1\otimes\lambda)$. By Lemma~\ref{lem:V}, the injective
homomorphism
\[
\mathrm{Ad}V\circ\widetilde{\pi^{\omega}}\rtimes(1\otimes\lambda):
\mathcal{A}^{\omega}_{\alpha}\rtimes_{\alpha^{\omega}}G\longrightarrow B((\mathcal{H}\otimes L^{2}(G))^{\omega})
\]
maps $af$ to
\[
(\tilde{\pi}\rtimes(1\otimes\lambda))^{\omega}(af)VV^{*}.
\]
Therefore, for each $af\in\mathcal{A}^{\omega}_{\alpha}\otimes_{\mathrm{alg}}C_{c}(G)$,
$\iota^{\omega}\rtimes u^{\omega}(af)VV^{*}=\mathrm{Ad}V\circ\widetilde{\pi^{\omega}}\rtimes(1\otimes\lambda)(af)$.
Thus, this equality holds for every $x\in\mathcal{A}^{\omega}_{\alpha}\rtimes_{\alpha^{\omega}}G$.
This proves that $\iota^{\omega}\rtimes u^{\omega}$ is injective.
\end{proof}

For abelian $G$, by a simple calculation, we obtain
\[
\iota^{\omega}\rtimes u^{\omega}\circ\widehat{\alpha^{\omega}}_{\gamma}([(a_{i}(n))]_{\omega}f_{i})
=\hat{\alpha}^{\omega}_{\gamma}\circ\iota^{\omega}\rtimes u^{\omega}([(a_{i}(n))]_{\omega}f_{i}),
\]
and hence, by density, $\iota^{\omega}\rtimes u^{\omega}$ is $\hat{G}$-equivariant. Since the dual
action is automatically continuous, in particular, we have
$\mathcal{A}^{\omega}_{\alpha}\rtimes_{\alpha^{\omega}}G\subset(\mathcal{A}\rtimes_{\alpha}G)^{\omega}_{\hat{\alpha}}$.
The question remains whether the reverse inclusion holds.

\section{Reverse inclusion}

The following diagram summarizes the relationships between ultraproducts of crossed products by
a discrete abelian group and the corresponding dual actions:
\[
\begin{tikzcd}
 (\mathcal{A}^{\omega}\rtimes_{\alpha^{\omega}}\Gamma)\rtimes_{\widehat{\alpha^{\omega}}}\hat{\Gamma}
    \arrow[r, equal]
    \arrow[d, hook]
 & \mathcal{A}^{\omega}\otimes\mathbb{K}
    \arrow[rd, "\sigma"]\\
 (\mathcal{A}\rtimes_{\alpha}\Gamma)^{\omega}_{\hat{\alpha}}\rtimes_{\hat{\alpha}^{\omega}}\hat{\Gamma}
    \arrow[r, hook]
 & (\mathcal{A}\rtimes_{\alpha}\Gamma\rtimes_{\hat{\alpha}}\hat{\Gamma})^{\omega}_{{\hat{\hat{\alpha}}}}
    \arrow[r, equal]
 & (\mathcal{A}\otimes\mathbb{K})^{\omega}_{\alpha\otimes \mathrm{Ad}\rho}
\end{tikzcd}
\]
A key step in the proof of Theorem~\ref{thm:Tomatsu}(2) is the equality
\[
(M \otimes B(L^{2}(G)))_{\alpha \otimes \operatorname{Ad} \rho}^{\omega} = M_{\alpha}^{\omega} \otimes B(L^{2}(G)),
\]
which follows from the tensor product decomposition for a type I factor $Q\subset M$, namely,
$M=Q\vee(Q^\prime \cap M)$ (cf.\ \cite[Lemma 2.8]{MT16}). However, in our situation, we do not
have sufficient information about the inverse of $\sigma$. Therefore, the von Neumann algebra
argument cannot be directly applied.

For a free ultrafilter $\omega$, the non-separability of $c_{0}(\mathbb{Z})^{\omega}$ can be proved
in the same way as for $\ell^{\infty}$.

\begin{lemma}\label{lem:c0}
The ultraproduct $c_0(\mathbb{Z})^\omega$ is not separable.
\end{lemma}

\begin{proof}
For $k\in\mathbb{N}$, let $\delta_k\in c_0(\mathbb{Z})$ denote the characteristic function of
$\{k\}$. For $I\subset\mathbb{N}$, define $e_I\in\ell^\infty(c_0(\mathbb{Z}))$ by
\[
e_I(n):=\sum_{k\in I,\ k\le n}\delta_k .
\]
 Let $I\neq J$
and put $k_0:=\min(I\triangle J)$.  Then
\[
\|e_I(n)-e_J(n)\|=
\begin{cases}
1 & (n\ge k_0),\\
0 & (n<k_0).
\end{cases}
\]
Since $\omega$ is free, $\{n: n\ge k_0\}\in\omega$. With
$q\colon\ell^\infty(c_0(\mathbb{Z}))\twoheadrightarrow c_0(\mathbb{Z})^\omega$ the quotient map,
this gives
\[
\|q(e_I)-q(e_J)\|=\lim_{n\to\omega}\|e_I(n)-e_J(n)\|=1 .
\]
Now let $S\subset c_0(\mathbb{Z})^\omega$ be dense. For each $I\in \mathcal{P}(\mathbb{N})$, choose
$s_I\in S$ with $\|q(e_I)-s_I\|<1/2$. If $I\neq J$, the triangle inequality gives
\[
\|s_I-s_J\|\ \ge\ \|q(e_I)-q(e_J)\|-\|q(e_I)-s_I\|-\|q(e_J)-s_J\|\ >\ 1-1/2-1/2=0 .
\]
Hence $I\mapsto s_I$ is injective, so $|S|\ge|\mathcal{P}(\mathbb{N})|=2^{\aleph_0}$. 
\end{proof}

It is natural to ask whether the reverse inclusion always holds. The main point is that for a
compact group, its dual is discrete, so the ultra-action is automatically continuous, which
enlarges the right-hand side.

\begin{proposition}\label{prop:counterexample}
There exist a locally compact abelian group $G$ and a C$^*$algebra $\mathcal{A}$ such that
\[
\mathcal{A}^{\omega}_{\alpha}\rtimes_{\alpha^{\omega}}G \neq (\mathcal{A}\rtimes_{\alpha}G)^{\omega}_{\hat{\alpha}}.
\]
\end{proposition}

\begin{proof}
Let us assume that $\mathcal{A}=\mathbb{C}$, $G=\mathbb{T}$, and let $\tau$ denote the trivial
action. By the Pontryagin duality theorem, we have
$C^{*}\mathbb{T}\simeq C_{0}(\mathbb{Z})=c_{0}(\mathbb{Z})$, which is separable. Moreover,
$C^{*}\mathbb{T}^{\omega}_{\tau}=C^{*}\mathbb{T}^{\omega}=c_{0}(\mathbb{Z})^{\omega}$. Hence, by
Lemma~\ref{lem:c0}, $\mathbb{C}\rtimes_{\tau}\mathbb{T}\subsetneq
(\mathbb{C}\rtimes_{\tau}\mathbb{T})^{\omega}_{\hat{\tau}}$.
\end{proof}

To avoid the situation in Proposition~\ref{prop:counterexample}, let us set that $G=\Gamma$ is a
discrete abelian group. In this case, does the reverse inclusion hold? For this purpose, we
introduce some notation and a few auxiliary lemmas.

\begin{notation}[{cf.\ \cite{Ped79}}]
For $s\in\Gamma$, set the closed subspace
\[
\mathcal{A}\rtimes_{\alpha}\Gamma(s):=\{x\in\mathcal{A}\rtimes_{\alpha}\Gamma:
\hat{\alpha}_{\gamma}(x)=\braket{s,\gamma}x \text{ for all } \gamma\}.
\]
For $x\in(\mathcal{A}\rtimes\Gamma)\setminus\{0\}$,
$F_{x}:\hat{\Gamma}\ni \gamma\mapsto\overline{\braket{s,\gamma}}\hat{\alpha}_{\gamma}(x)$ is
continuous. Since $\mathrm{supp}(F_{x})=\hat{\Gamma}$ is compact, the Haar integral of $F_{x}$
exists. More precisely, by Lemma~\ref{lem:integral}, we can define
\[
E_{s}(x):=\int_{\hat{\Gamma}}\overline{\braket{s,\gamma}}\hat{\alpha}_{\gamma}(x)\,d\gamma
\in\mathcal{A}\rtimes\Gamma.
\]
Moreover, it is straightforward to check that $E_{s}$ is a map onto
$\mathcal{A}\rtimes_{\alpha}\Gamma(s)$ and satisfies $\|E_{s}(x)\|\leq\|x\|$ and $E_{s}^{2}=E_{s}$.
\end{notation}

\begin{lemma}\label{lem:spectral}
$\mathcal{A}\lambda_{s}=\mathcal{A}\rtimes\Gamma(s)$ for $s\in\Gamma$.
\end{lemma}

\begin{proof}
Take $x\in\mathcal{A}\rtimes\Gamma(s)$. Then
$\hat{\alpha}_{\gamma}(x\lambda^{*}_{s})=\braket{s,\gamma}x\overline{\braket{s,\gamma}}\lambda^{*}_{s}
=x\lambda^{*}_{s}$. Thus,
$x\lambda^{*}_{s}\in\mathcal{A}\rtimes_{\alpha}\Gamma(e)=(\mathcal{A}\rtimes_{\alpha}\Gamma)^{\hat{\alpha}}
=\mathcal{A}$. Therefore, $x\in\mathcal{A}\lambda_{s}$.
\end{proof}

By the preceding lemma, we know that the C$^*$subalgebra
$\overline{\mathrm{span}}_{s}\mathcal{A}\rtimes_{\alpha}\Gamma(s)\subset \mathcal{A}\rtimes_{\alpha}\Gamma$
actually coincides with the entire crossed product $\mathcal{A}\rtimes_{\alpha}\Gamma$.

\begin{lemma}\label{lem:span}
For the compact group action
$\hat{\alpha}^{\omega}:\hat{\Gamma}\curvearrowright(\mathcal{A}\rtimes_{\alpha}\Gamma)^{\omega}_{\hat{\alpha}}$,
we have
$(\mathcal{A}\rtimes_{\alpha}\Gamma)^{\omega}_{\hat{\alpha}}=\overline{\mathrm{span}}_{s}(\mathcal{A}\rtimes\Gamma)^{\omega}(s)$.
\end{lemma}

\begin{proof}
The continuity of $\hat{\alpha}^{\omega}$ allows us to consider
$E_{s,\omega}:(\mathcal{A}\rtimes_{\alpha}\Gamma)^{\omega}_{\hat{\alpha}}\twoheadrightarrow(\mathcal{A}\rtimes\Gamma)^{\omega}(s)$;
$E_{s,\omega}(x)=\int_{\hat{\Gamma}}\overline{\braket{s,\gamma}}\hat{\alpha}^{\omega}_{\gamma}(x)\,d\gamma$,
which is well defined by Lemma~\ref{lem:integral}.

Let us assume that
$(\mathcal{A}\rtimes_{\alpha}\Gamma)^{\omega}_{\hat{\alpha}}\setminus\overline{\mathrm{span}}_{s}(\mathcal{A}\rtimes\Gamma)^{\omega}(s)$
is nonempty, and that $x$ belongs to it. By the Hahn--Banach theorem, there exists a linear
functional $\varphi\in((\mathcal{A}\rtimes_{\alpha}\Gamma)^{\omega}_{\hat{\alpha}})^{*}$ such that
$\varphi(x)\neq 0$ and $\varphi(y)=0$ for all $y\in\mathrm{span}_{s}(\mathcal{A}\rtimes\Gamma)^{\omega}(s)$.
Hence, for all $s\in\Gamma$,
$\varphi(E_{s,\omega}(x))=\int_{\hat{\Gamma}}\overline{\braket{s,\gamma}}\varphi(\hat{\alpha}^{\omega}_{\gamma}(x))\,d\gamma=0$.
By the uniqueness of Fourier coefficients, $\varphi(\hat{\alpha}^{\omega}_{\gamma}(x))=0$ for all
$\gamma$. In particular, taking $\gamma=e$ gives a contradiction.
\end{proof}

Note that $(x_{n})_{n}\in \ell^{\infty}(\mathcal{A}\rtimes\Gamma)$ is uniformly bounded. Since
$E_{s}$ is bounded, we can define
$E^{\omega}_{s}:(\mathcal{A}\rtimes\Gamma)^{\omega}_{\hat{\alpha}}\to\mathcal{A}\rtimes\Gamma(s)^{\omega}$;
$E^{\omega}_{s}([(x_{n})_{n}]_{\omega}):=[(E_{s}(x_{n}))_{n}]_{\omega}$. Furthermore, by
Lemma~\ref{lem:spectral}, the image of $E_{s}^{\omega}$ is contained in
$(\mathcal{A}\lambda_{s})^{\omega}=\mathcal{A}^{\omega}\lambda_{s}$.

\begin{theorem}\label{thm:reverse}
For a discrete abelian group $\Gamma$, we have
$\mathcal{A}^{\omega}\rtimes_{\alpha^{\omega}}\Gamma=(\mathcal{A}\rtimes_{\alpha}\Gamma)^{\omega}_{\hat{\alpha}}$.
\end{theorem}

\begin{proof}
For $x\in(\mathcal{A}\rtimes\Gamma)^{\omega}(s)$, there exists
$[(y_{n})_{n}]_{\omega}\in (\mathcal{A}\rtimes\Gamma)^{\omega}_{\hat{\alpha}}$, where the lift
$(y_{n})$ of $x$ lies in the continuous part of $\ell^{\infty}(\mathcal{A}\rtimes\Gamma)$, such
that $E_{s,\omega}([(y_{n})_{n}]_{\omega})=x$. Then
\[
x=\int_{\hat{\Gamma}}\overline{\braket{s,\gamma}}\hat{\alpha}^{\omega}_{\gamma}([(y_{n})_{n}]_{\omega})\,d\gamma
=\int[(\overline{\braket{s,\gamma}}\hat{\alpha}_{\gamma}(y_{n}))_{n}]_{\omega}\,d\gamma.
\]
Since the restriction of the quotient map
$[\cdot]_{\omega}:\ell^{\infty}(\mathcal{A}\rtimes\Gamma)\twoheadrightarrow (\mathcal{A}\rtimes\Gamma)^{\omega}$
to the continuous part is also bounded (or equivalently, by using Lemma~\ref{lem:integral} on the
properties of integrals of C$^*$algebra-valued functions), the above formula equals
\[
\Bigl[\int(\overline{\braket{s,\gamma}}\hat{\alpha}_{\gamma}(y_{n}))_{n}\,d\gamma\Bigr]_{\omega}
=E_{s}^{\omega}([(y_{n})_{n}]_{\omega})\in\mathcal{A}^{\omega}\lambda_{s}.
\]
This shows that
$x\in\overline{\mathrm{span}}_{s}\mathcal{A}^{\omega}\lambda_{s}=\mathcal{A}^{\omega}\rtimes\Gamma$.
Hence, by Lemma~\ref{lem:span}, we obtain the desired conclusion.
\end{proof}

\subsection*{Acknowledgements}
The author is grateful to Toshihiko Masuda for his support and valuable comments on this paper as
my supervisor. The author would also like to thank Miho Mukohara for helpful discussions
regarding this work.

\bigskip
\noindent Graduate School of Mathematics, Kyushu University, Motooka 744, Nishi-ku, Fukuoka
819-0395, Japan\\
\noindent Email address: \texttt{canenblluba@gmail.com}


\begin{thebibliography}{99}
\bibitem[AH14]{AH14}
H.~Ando and U.~Haagerup, Ultraproducts of von Neumann algebras, J. Funct. Anal. \textbf{266}
(2014), 6842--6913.
\bibitem[BO08]{BO08}
N.~P.~Brown and N.~Ozawa, C$^*$-algebras and finite-dimensional approximations, Graduate Studies
in Mathematics, vol.~88, American Mathematical Society, Providence, RI, 2008.
\bibitem[DW07]{DW07}
D.~Williams, Crossed products of C$^*$-algebras, Mathematical Surveys and Monographs, vol.~134,
American Mathematical Society, Providence, RI, 2007.
\bibitem[Kis96]{Kis96}
A.~Kishimoto, A Rohlin property for one-parameter automorphism groups, Comm. Math. Phys.
\textbf{179} (1996), no.~3, 599--622.
\bibitem[MT16]{MT16}
T.~Masuda and R.~Tomatsu, Rohlin flows on von Neumann algebras, Mem. Amer. Math. Soc.
\textbf{244} (2016), no.~1153, ix+111 pp.
\bibitem[Ped79]{Ped79}
G.~K.~Pedersen, C$^*$-algebras and their automorphism groups, London Mathematical Society
Monographs, vol.~14, Academic Press, London, 1979.
\bibitem[Pis09]{Pis09}
G.~Pisier, On the lifting property for C$^*$-algebras, J. Noncommut. Geom. \textbf{16} (2022),
no.~3, 967--1006.
\bibitem[Rae88]{Rae88}
I.~Raeburn, On crossed products and Takai duality, Proc. Edinb. Math. Soc. \textbf{31} (1988),
321--330.
\bibitem[Tom17]{Tom17}
R.~Tomatsu, Ultraproducts of crossed product von Neumann algebras, Illinois J. Math. \textbf{61}
(2017), no.~3--4, 275--286.
\end{thebibliography}
\end{document}